\documentclass[a4paper,twoside,11pt]{book}

\usepackage[latin1]{inputenc}
\usepackage[T1]{fontenc}
\usepackage[english,french]{babel}

\usepackage{graphicx}

\usepackage{amsmath,amsthm,amsfonts,amssymb,amsxtra}

\usepackage{scrtime}

\usepackage{fancyhdr}
\newcommand{\styledehead}{\small\em}



\fancypagestyle{plain}{%
\fancyhf{} 
}

\newcommand{\beq}{\begin{equation}}
\newcommand{\eeq}{\end{equation}}
\newcommand{\bea}{\begin{eqnarray}}
\newcommand{\eea}{\end{eqnarray}}
\newcommand{\beas}{\begin{eqnarray*}}
\newcommand{\eeas}{\end{eqnarray*}}

\newtheorem{theorem}{Theorem}[section]

\newtheorem{corollary}[theorem]{Corollary}

\newtheorem{theoreme-fr}[theorem]{Théorème}
\newtheorem{definition-fr}[theorem]{Définition}
\newtheorem{proposition-fr}[theorem]{Proposition}
\newtheorem{corollaire-fr}[theorem]{Corollaire}
\newtheorem{lemme-fr}[theorem]{Lemme}
\newtheorem{remarque-fr}[theorem]{Remarque}
\newtheorem{exemple-fr}[theorem]{Exemple}
\newtheorem{exemples-fr}[theorem]{Exemples}

\usepackage{epsf}
\usepackage{bm}
\usepackage{makeidx}
\usepackage{latexsym}
\usepackage{amssymb}
\usepackage[latin1]{inputenc}
\usepackage{amsfonts}
\usepackage{amsmath}
\usepackage{amssymb}

\RequirePackage{amsopn} \RequirePackage{amsfonts}
\RequirePackage{amsthm}

\renewcommand{\thesection}{\arabic{section}}

\renewcommand{\thetheorem}{\thesection.\arabic{theorem}}

\begin{document}
\chapter*{Large deviations for directed percolation on a
thin rectangle.}

JEAN-PAUL IBRAHIM \\
\textit{University of Toulouse, France}

\newcounter{remarque}
\addtocounter{remarque}{1}


\

\

\
\begin{small}
Following the recent investigations of Baik and Suidan in
\cite{baik2005gcl} and Bodineau and Martin  in
\cite{bodineau2005upl}, we prove large deviation properties for a
last-passage percolation model in $\mathbb{Z}^{2}_{+}$ whose paths
are close to the axis. The results are mainly obtained when the random
weights are Gaussian or have a finite moment-generating function and rely,
as in \cite{baik2005gcl} and
\cite{bodineau2005upl}, on an embedding in Brownian paths and the
KMT approximation. The study of the subexponential
case completes the exposition.
\end{small}

\section*{Introduction}

Random matrix theory has developed extensively in the last several
decades following the pioneering results by E. Wigner in the
fifties. Gaussian models attracted a lot of attention, among them
the Gaussian Unitary Ensemble (GUE). In this example, the knowledge
of the joint distribution of the eigenvalues allowed for a rather
complete understanding of both their global and local behaviors. In
particular, the limiting behavior of the largest eigenvalue gave
rise to the famous Tracy-Widom distribution
 \cite{joh2,mehta,tw}. Several random growth models,
such as the longest increasing subsequence and the corner growth
models, have been shown to develop a similar behavior relying on a
common determinantal structure \cite{joh1,konig}. In particular, the
last-passage percolation or so-called corner growth model (see below
for a precise description) has been deeply studied by Johansson in
\cite{joh1}. For geometric or exponential random variables (the only
cases leading to a determinantal description), Johansson established
both fluctuations and large deviation asymptotics similar to the
ones for the GUE random matrix model. Following the recent
investigations by Baik and Suidan \cite{baik2005gcl} and Bodineau
and Martin \cite{bodineau2005upl} at the level of fluctuations, the
present paper deals with large deviations for the random growth
model for more general random variables but on rectangles such that
one side is asymptotically negligible with respect to the other at a
given rate. The main results concern Gaussian random variables and
random weights having a finite moment-generating function. Somewhat
surprisingly, the rate may be shown to be larger than the one for
the fluctuations. The comparison method used in this work is
basically inspired from \cite{baik2005gcl} and
\cite{bodineau2005upl} and relies similarly on an embedding in
Brownian paths.

Recall first the basic corner growth model under study. It can be
described as directed paths in the lattice $\mathbb{Z}^{2}_{+}$
going from $(1,1)$ to  $(N,k)\in\mathbb{Z}^{2}_{+}$ where only up
and right steps are allowed. More precisely, denoting by $\Pi(N,k)$
the set of all such paths, a path $\pi\in\Pi(N,k)$ is called an
\textit{up/right} path and is defined as a collection of sites
$\left\{(i_{l},j_{l})\right\}_{l=1}^{N+k-1}$ satisfying
$(i_{1},j_{1})=(1,1)$,
 $(i_{N+k-1},j_{N+k-1})=(N,k)$ and $(i_{l+1},j_{l+1})-(i_{l},j_{l})$ is either $(1,0)$ or $(0,1)$.
The main random variable under consideration is the last-passage time defined by
\begin{equation*}
G(N,k)=\underset{\pi\in\Pi(N,k)}{\max}\bigg
\{\sum_{(i,j)\in\pi}X_{i}^{(j)} \bigg \}
\end{equation*}
where the $X_{i}^{(j)}$'s are $i.i.d.$ random variables. As an
alternate description, set $U(N,k)$ as the subset of
$\mathbb{R}_{+}^{k+1}$ given by
$$
U(N,k) = \big \{ {u=(u_{0},u_{1},\ldots,u_{k})\in
\mathbb{R}_{+}^{k+1}; 0=u_{0}\leq u_{1}\leq \ldots \leq u_{k}=N}
\big \} .
$$
Then
$$
G(N,k)=\underset{u\in U(N,k)}{\sup}\bigg\{ \sum_{r=1}^{k} \big  [
S_{\lfloor u_{r}\rfloor}^{(r)}-S_{\lfloor
u_{r-1}\rfloor-1}^{(r)}\big ]\bigg \}
$$
where $S_{m}^{(r)}=\sum_{i=1}^{m}X^{(r)}_{i}$ with the convention
$S_{-1}^{(r)}=S_{0}^{(r)}=0$. This follows from the fact that
$U(N,k)\equiv\Pi(N,k)$ when $u \in \mathbb{Z}_{+}^{k+1}$. Actually,
every $u\in U(N,k)\cap \mathbb{Z}_{+}^{k+1}$ maps to a unique path
$\pi \in \Pi(N,k)$ whose $i^{th}$ up-jump occurs on $u_{i}$. On the
other hand, each path $\pi$ is characterized by its up-step sites.
It will be more appropriate
 to adopt the second form of $G(N,k)$ in order to compare it later
with the Brownian last-passage percolation model.


\section{Results}\label{sec1}
Having introduced the model, recall briefly some of its properties
established by Johansson in \cite{joh1} in the particular case of
geometric or exponential distributions. The key of the results in
this case relies on the explicit description
 of the last-passage time distribution $G(N,k)$. When weights are $i.i.d.$ geometric
 with parameter $q\in(0,1)$, we have
\begin{equation*}
\mathbb{P}[G(N,k)\leq t]=\frac{1}{Z_{N,k}}\mathop{\sum_{
h\in\mathbb{N}^{k}}}_{\max\{h_{i}\}\leq t+k-1}\prod_{1\leq i<j\leq
k}(h_{i}-h_{j})^{2}\prod_{i=1}^{k}\binom{h_{i}+N-k}{h_{i}}q^{h_{i}},
\end{equation*}
where $N\geq k$ and $Z_{N,k}$ is the normalizing constant. Using
results from logarithmic potential theory, Johansson described in
\cite{joh1} the large deviation behaviors of $G(N,k) $ when $ k \sim
N$. Namely, he obtained that, for $\gamma\geq1$, there exist two
functions $i(\varepsilon)$ and $l(\varepsilon)$ such that for any
$\varepsilon>0$,
$$
\lim_{N\rightarrow\infty}  \frac{1}{N} \log \mathbb{P}[G([\gamma
N],N)\geq N(\omega(\gamma,q)+\varepsilon)]=-i(\varepsilon)
$$
and
$$
\lim_{N\rightarrow\infty}  \frac{1}{N^{2}} \log \mathbb{P}[G([\gamma
N],N)\leq N(\omega(\gamma,q)-\varepsilon)]=-l(\varepsilon).
$$
The functions $l(x)$ and $i(x)$ are positive for every $x>0$.
Furthermore,
$$
\lim_{N\rightarrow\infty} \frac{1}{N} \, \mathbb{E}[G([\gamma
N],N)]=\omega(\gamma,q):=\frac{(1+\sqrt{\gamma q})^{2}}{1-q}-1.
$$

Using the asymptotics of the Meixner orthogonal polynomial ensemble,
Johansson \cite{joh1} further established the fluctuations of
$G(N,k)$ at the Tracy-Widom GUE rate. He proved that for
$\gamma\geq1$ and $s\in\mathbb{R}$,
$$
\lim_{N\rightarrow\infty}\mathbb{P}\Big[\frac{G([\gamma
N],N)-N\omega(\gamma,q)}{\sigma(\gamma,q)N^{1/3}}\leq
s\Big]=F_{\textrm{TW}}(s),
$$
where
$$
\sigma(\gamma,q)=\frac{q^{1/6}\gamma^{-1/6}}{1-q}(\sqrt{\gamma}+\sqrt{q})^{2/3}(1+\sqrt{\gamma
q})^{2/3}
$$
and $F_{\textrm{TW}}(s)$ is the distribution function of the
Tracy-Widom law (see \cite{tw}). Replacing geometric weights with
exponential ones gives similar results since an exponential
distribution can be seen as the limit of a rescaled geometric one.
See \cite{joh1} for the precise formulas.

\

Recently, Bodineau and Martin \cite{bodineau2005upl} and Baik and
Suidan \cite{baik2005gcl} studied the same model when paths are
close to the axis, \textit{i.e.} $k=o (N^{\alpha})$ for some $\alpha
< 1$ but allowing more general distributions. The authors used a
coupling with the Brownian trajectories through the following
Brownian last-passage percolation. Letting $(B_{t}^{(r)})_{r\geq 1}$
be a sequence of independent Brownian motions, set
$$
L(N,k)=\underset{u\in U(N,k)}{\sup}\bigg\{ \sum_{r=1}^{k}\big [
B_{u_{r}}^{(r)}-B_{u_{r-1}}^{(r)}\big ]\bigg \}.
$$
It has been proved in  \cite{baryshnikov,gravner,connell_yor} that
$L(1,k)$ has the same distribution as the largest eigenvalue of a
$k\times k$ rescaled GUE random matrix. As a consequence of the
fluctuation result for the GUE model, it follows that
$$
k^{1/6} \big [L(1,k)-2\sqrt{k} \big
]\overset{d}{\longrightarrow}F_{\textrm{TW}}.
$$
Using this result and a comparison between the continuous model with
Brownian paths and the discrete one with random weights, the authors
of \cite{baik2005gcl} and  \cite{bodineau2005upl} deduced
fluctuation properties of the corner growth model for rather general
random variables. However, the embedding in the Brownian paths
requires to restrict the paths on small rectangles. For example, in
\cite{bodineau2005upl}, the discrete and the continuous models were
coupled using the Koml\'os-Major-Tusn\'ady (KMT) approximation which
couples random walks with Brownian motion. The authors proved that
if the weights satisfy $\mathbb{E}|X_{i}^{(j)}|^{p}<\infty$ for some
$p>2$, setting $\mu=\mathbb{E}X_{i}^{(j)}$ and
$\sigma^{2}=\textrm{var}(X_{i}^{(j)})$, then for all
$\alpha<\frac{6}{7}(\frac{1}{2}-\frac{1}{p})$,
\begin{equation*}
\frac{G(N,\lfloor N^{\alpha}\rfloor)-N\mu-2\sigma
N^{\frac{1+\alpha}{2}}}{\sigma
N^{\frac{1}{2}-\frac{\alpha}{6}}}\overset{d}{\longrightarrow}F_{\textrm{TW}}.
\end{equation*}
If the random variables $X_{i}^{(j)}$ have all moments,
\textit{i.e.} $p=\infty$, then $\alpha$ is lower than $3/7$. This is
true when the weights are Gaussian or are bounded for example. It is
not known how optimal this rate could be: the authors in
\cite{bodineau2005upl} think that such a result might hold, for some
independence reasons, when $\alpha<3/4$. However, they do not give a
complete proof. In \cite{baik2005gcl}, the authors compared the
discrete and continuous model via the Skorokhod embedding theorem in
order to obtain almost the same results. Lately, Suidan in
\cite{suidan2006rtc} produced another proof of the last theorem when
the variables have a third moment. He compared two discrete directed
percolation models using a theorem of Chatterjee
\cite{chatterjee2005sit}. The fluctuation properties for the first
(with geometric distribution) lead him to similar ones for the
second.

\

In this paper, we follow the comparison methods of
\cite{bodineau2005upl} and \cite{baik2005gcl} to establish large
deviations limit theorems for directed percolation models on thin
rectangles mainly for Gaussian distribution. An extension for random
weights with finite moment-generating function aroud zero will be
also proved. We rely similarly on the corresponding results for the
Brownian percolation model. Namely, as a consequence of the GUE
random matrix interpretation \cite{gravner,connell_yor}, for all
$\varepsilon>0$,
\begin{equation}\label{E2}
 \lim_{k \to \infty} \frac{1}{k} \, \log{\mathbb{P}\big [L(1,k) \geq
2\sqrt{k}(1+\varepsilon)\big]}= -J_{GUE}(\varepsilon)
\end{equation}
and
\begin{equation}\label{E3}
 \lim_{k \to \infty}\frac{1}{k^{2}}\, \log{\mathbb{P}\big [L(1,k) \leq
2\sqrt{k}(1-\varepsilon)\big ]}= -I_{GUE}(\varepsilon).
\end{equation}
The two functions $J_{GUE}(x)$ and $I_{GUE}(x)$ are both positive
for every positive $x$.  $J_{GUE}$ can be computed explicitly (see
\cite{benarous}) as
$$
J_{GUE}(\varepsilon)=4\displaystyle\int_{0}^{\varepsilon}\sqrt{x(x+2)}~dx.
$$
To the best of our knowledge, there is no explicit form for
$I_{GUE}$. This function appears in the logarithmic potential theory
and it represents physically the minimal potential energy  of
charges on a one-dimension conductor exposed to an external field
(see \cite{saff}). In this work, we do not need the explicit form of
$I_{GUE}$. However, its continuity, proved at the end of Section
\ref{sec2}, will be necessary for the proof.

The following three theorems are the main results of this paper.
Despite some similarity in their proofs, the second one requires
more work. One surprising feature is that the rate $\alpha $ in the
Gaussian case is less then one (compared to $\alpha < 3/7$ for the
fluctuation result).


Throughout the article, $k$ and $N$ are two integers which depend on each other,
\textit{i.e.} $N=N(k)=N_k$. We assume that $k=\mathrm{o}(N_k)$ and we let $k$ goes
to infinity. For simplicity, we write $N$ instead of
$N_k$ throughout the proofs.

\begin{theorem}\label{th1}
Assume that the variables $(X_{i}^{(j)})_{i,j=1}^{\infty}$ are
i.i.d. standard normal random variables. Assume further that $k=\mathrm{o}
(\frac{N_k}{\log{N_k}})$. Then, for all $\varepsilon
>0$,
\begin{equation*}
 \lim_{k\to \infty} \frac{1}{k} \, \log{\mathbb{P}\big [G(N_k,k)
\geq 2\sqrt{N_kk}(1+\varepsilon)\big]}= -J_{GUE}(\varepsilon).
\end{equation*}
\noindent On the left of the mean, we have for
$k=\mathrm{o}(N_k^{\frac{1}{2}})$,
\begin{equation*}
 \lim_{k\to \infty} \frac{1}{k^{2}} \, \log{\mathbb{P}\big [G(N_k,k) \leq
2\sqrt{N_kk}(1-\varepsilon)\big]}= -I_{GUE}(\varepsilon).
\end{equation*}
\end{theorem}
In the second statement, we replace Gaussian variables with weights
having finite exponential moments. Loosing the Gaussian assumption
will complicate the coupling and reduce the size of the rectangles.
We denote by $X$ a random variable having the common law of the
$i.i.d.$ variables in the sequence $(X_{i}^{(j)})_{i,j=1}^{\infty}$.


\begin{theorem}\label{th2}
Assume that the variables $(X_{i}^{(j)})_{i,j=1}^{\infty}$ are
i.i.d. random variables such that $\mathbb{E}X=0$ and
$\mathbb{E}X^{2}=1$. Assume further that there exit $\mu_0>0$ such that for all
$\mu<\mu_0$,
\begin{equation}\label{E46bis}
\mathbb{E}\exp{(\mu|X|)}<+\infty.
\end{equation}. If $k=\mathrm{o}(\frac{N_k}{(\log{N_k})^{2}})$, then for all $\varepsilon
>0$,
\begin{equation*}
 \lim_{k\to \infty} \frac{1}{k} \, \log{\mathbb{P}\big [G(N_k,k)
\geq 2\sqrt{N_kk}(1+\varepsilon)\big]}= -J_{GUE}(\varepsilon).
\end{equation*}
\noindent Similarly, if $k=\mathrm{o}(N_k^{\frac{1}{3}})$,  for all $\varepsilon>0$,
\begin{equation*}
 \lim_{k\to \infty} \, \frac{1}{k^{2}}~\log{\mathbb{P}\big [G(N_k,k) \leq
2\sqrt{N_kk}(1-\varepsilon)\big]}= -I_{GUE}(\varepsilon) .
\end{equation*}
\end{theorem}
The proof of Theorem \ref{th2} relies on the
Koml\'os-Major-Tusn\'ady approximation for the sums of $i.i.d$
centered random variables with finite exponential moments, see
\cite{komlos}. The following theorem deals with a particular class
of subexponential weights. We make use of the Skorokhod embedding
theorem to obtain the following results.


\begin{theorem}\label{th7}
Assume that the variables $(X_{i}^{(j)})_{i,j=1}^{\infty}$ are
i.i.d. random variables satisfying $\mathbb{E}X=0$ and
$\mathbb{E}X^{2}=1$. Furthermore, assume that there exit $\mu>0$ and
$0<\gamma<1$ such that
\begin{equation}\label{E46}
\mathbb{E}\exp{(\mu|X|^{\gamma})}<+\infty.
\end{equation}
 If $k=\mathrm{o}(N_k^{\alpha})$ with $\alpha< \frac{\gamma}{2\gamma+2}$, then, for
all $\varepsilon
>0$,
\begin{equation*}
 \lim_{k\to \infty} \frac{1}{k} \, \log{\mathbb{P}\big [G(N_k,k)
\geq 2\sqrt{N_kk}(1+\varepsilon)\big]}= -J_{GUE}(\varepsilon).
\end{equation*}
\noindent Similarly, if $k=\mathrm{o}(N_k^{\alpha})$ with $\alpha<\frac{\gamma}{5\gamma+4}$,  for all
$\varepsilon>0$,
\begin{equation*}
 \lim_{k\to \infty} \, \frac{1}{k^{2}}~\log{\mathbb{P}\big [G(N_k,k) \leq
2\sqrt{N_kk}(1-\varepsilon)\big]}= -I_{GUE}(\varepsilon) .
\end{equation*}
\end{theorem}
The results in Theorem \ref{th7} cover in particular the examples of
Weibull and L\'evy distributions. Notice that in Theorem \ref{th7}
we can take $\gamma\geq1$. However, the result is worthless because
of Theorem \ref{th2}. Actually, the KMT approximation is more
efficient than the Skorokhod embedding theorem for $i.i.d.$ random
variables with finite exponential moments.

\

The asymptotic results above hold also for geometric and exponential
weights. Moreover, the rectangle width in this case, has to be only
negligible with respect to its length  since the rate functions of the
Laguerre ensemble converge to the GUE ones when $k=\mathrm{o}(N_k)$.
The reader can see \cite{joh1} and \cite{ledoux} for rigorous
results.

 Non-asymptotic bounds for the preceding models can be deduced
from the previous theorems proofs. The rectangle width for small deviations
matches in this case the fluctuation results. For this we use analogous
deviation inequalities to the right of the mean obtained for the largest
eigenvalue of the GUE , see
\cite{ledoux}. To the left of the mean, we use recent deviation results
for the largest eigenvalue of the GUE obtained by Ledoux and Rider \cite{ledouxrider}.


\begin{theorem}\label{th5}
Assume that the variables $(X_{i}^{(j)})_{i,j=1}^{\infty}$ are
i.i.d. standard normal random variables, and that
 $k=N_k^{\alpha}$ with $\alpha< \frac{3}{7}$. Then, there exists a
positive constant $C_{\alpha}$  depending only on $\alpha$ such
that, for all $0<\varepsilon <1$,
\begin{equation}
\mathbb{P}\big [G(N,k) \geq 2\sqrt{N k}(1+\varepsilon)\big]\leq
C_{\alpha}\exp{\Big (-
\frac{k\varepsilon^{3/2}}{C_{\alpha}}\Big)}.
\end{equation}
On the left of the mean, we have
\begin{equation}
\mathbb{P}\big [G(N,k) \leq 2\sqrt{N k}(1-\varepsilon)\big]\leq
C_{\alpha}\exp{\Big (-
\frac{k^2\varepsilon^{3}}{C_{\alpha}}\Big)}.
\end{equation}
\end{theorem}

\

\noindent For random weights with finite moment-generating function,
we have similar results.

\begin{theorem}\label{th6}
Assume that the variables $(X_{i}^{(j)})_{i,j=1}^{\infty}$ are
i.i.d. random variables such that $\mathbb{E}X=0$ and
$\mathbb{E}X^{2}=1$. Assume further that there exit $\mu_0>0$ such that for all
$\mu<\mu_0$,
\begin{equation}
\mathbb{E}\exp{(\mu|X|)}<+\infty.
\end{equation}. If
 $k=N_k^{\alpha}$ with $\alpha< \frac{3}{7}$, then there exists a
positive constant $C_{\alpha}$  depending $\alpha$ and the
distribution of $X$ such that, for all $0<\varepsilon <1$,
\begin{equation}\label{E8}
\mathbb{P}\big [G(N_k,k) \geq 2\sqrt{N_kk}(1+\varepsilon)\big]\leq
C_{\alpha}\exp{\Big (- \frac{k\varepsilon^{3/2}}{C_{\alpha}}\Big)}.
\end{equation}
Similarly, if
 $k=N_k^{\alpha}$ with $\alpha< \frac{1}{3}$, then
 \begin{equation}\label{E8biss}
\mathbb{P}\big [G(N_k,k) \leq 2\sqrt{N_kk}(1-\varepsilon)\big]\leq
C_{\alpha}\exp{\Big (- \frac{k^2\varepsilon^{3}}{C_{\alpha}}\Big)}.
\end{equation}
\end{theorem}
Large deviation inequalities for large $\varepsilon>1$ hold for the
optimal rate $\alpha=1$ using some concentration arguments. We refer
to \cite{ledoux} for the results and the proofs.

As the reader can notice, if $X$ satisfies the condition
(\ref{E46}), the same exponential inequalities as in Theorem
\ref{th6} can be obtained. In this case, the positive constant will
depend on $\alpha$, $\mu$ and $\gamma$. A smaller $\alpha$ will be
also necessary. The precise calculations are left to the reader.

Let us notice here that if the weights have a finite $p$-th moment with $p\geq2$,
we still obtain the same deviation results but for $k=\mathrm{o}(\log{(N^{p-1})})$.
In this case, we use again the KMT approximation for random weights with finite $p$-th moment, see \cite{komlos,sakhanenko}.
\

We strongly believe that the preceding results hold on a wider
rectangle. However, the method used here does not allow us to
improve the rectangle width. Theorem \ref{th1} will be proved in
Section \ref{sec2} while Theorem \ref{th2} will be proved in Section
\ref{sec3}. The case of the particular subexponential weights given
in Theorem \ref{th7} will be discussed in section \ref{sec4}.
Theorem \ref{th5} and Theorem \ref{th6} will be addressed in Section
\ref{sec5} on the basis of the preceding results and proofs.


\section{Proof of Theorem \ref{th1}.}\label{sec2}
\setcounter{equation}{0} We assume throughout the rest of the paper
and without loss of generality that $k$ and $N$ are two positive integers
and we write $N$ instead of $N_k$ and $N(k)$. As claimed
before, to prove Theorem \ref{th1}, we compare $G(N,k)$ and
$L(N,k)$. To do so, let for any $\varepsilon>0$,
$$
A=\big \{G(N,k)\geq 2\sqrt{N k}(1+\varepsilon)\big \}
$$
and
$$
B=\big \{\big|G(N,k)-L(N,k)\big|\geq
2\sqrt{N k}(\varepsilon-\varepsilon_{1})\big \}
$$
where $0<\varepsilon_{1}<\varepsilon$. Clearly,
\begin{equation}\label{E10}
\mathbb{P}[A]\leq \mathbb{P}\big [L(N,k)\geq
2\sqrt{N k}(1+\varepsilon_{1})\big ]+\mathbb{P}[B]
\end{equation}
and
\begin{equation}\label{E11}
\mathbb{P}[A]\geq \mathbb{P}\big [L(N,k)\geq
2\sqrt{N k}(1+2\varepsilon-\varepsilon_{1})\big ]-\mathbb{P}[B].
\end{equation}
Moreover, for every $\eta >0$,
\begin{equation}\label{E11bis}
\mathbb{P}\big [L(N,k)\geq 2\sqrt{N k}(1+\eta)\big
]=\mathbb{P}\big[L(1,k)\geq 2\sqrt{k}(1+\eta)\big]
\end{equation}
as a consequence of the Brownian scaling
$\sqrt{N}L(1,k)\overset{d}{=}L(N,k)$. To evaluate $\mathbb{P}[B]$,
we couple $G(N,k)$ and $L(N,k)$ by letting
$X_{i}^{(j)}=B^{(j)}_{i}-B^{(j)}_{i-1}$ for all $i,j\geq 1$ so that
the sequence $(X_{i}^{(j)})_{i,j=1}^{\infty}$ is $i.i.d.$ with
standard normal distribution. When comparing $G(N,k)$ and $L(N,k)$,
it is obvious that most of the variables will vanish. More
precisely, repeating the computation done by Bodineau and
 Martin in Section 2 of \cite{bodineau2005upl}, we get, by
 letting $B^{(r)}_{-1}=0$,
\begin{eqnarray}\label{E12}
&&\big|G(N,k)-L(N,k)\big|\nonumber\\ \nonumber\\
&&=\bigg|\displaystyle\sup_{u\in U(N,k)}\sum_{r=1}^{k}\Big [
S_{\lfloor u_{r}\rfloor}^{(r)}-S_{\lfloor
u_{r-1}\rfloor-1}^{(r)}\Big ]-\underset{u^{'}\in U(N,k)}{\sup}
\sum_{r=1}^{k}\Big[B_{u^{'}_{r}}^{(r)}-B_{u^{'}_{r-1}}^{(r)}\Big]\bigg|\nonumber\\ \nonumber\\
\nonumber\\
&&\leq \displaystyle\sup_{u\in U(N,k)}
\sum_{r=1}^{k}\Big[\Big|S_{\lfloor u_{r}\rfloor}^{(r)}-B_{\lfloor
u_{r}\rfloor}^{(r)}\Big|+\big|S_{\lfloor
u_{r-1}\rfloor-1}^{(r)}-B_{\lfloor
u_{r-1}\rfloor-1}^{(r)}\Big|\nonumber\\&&
\,\,\,\,\,\,\,\,\,\,\,\,\,\,\,\,\,\,\,\,\,\,\,\,\,\,\,\,\,\,\,\,+\Big|B_{\lfloor
u_{r}\rfloor}^{(r)}-B_{u_{r}}^{(r)}\Big|+\Big|B_{\lfloor
u_{r-1}\rfloor-1}^{(r)}-B_{u_{r-1}}^{(r)}\Big|\Big]\nonumber \\
\nonumber\\ \nonumber\\
&&\leq 2\displaystyle\sum_{r=1}^{k}
\bigg(\displaystyle\max_{i=1,\ldots,N}\big|S_{i}^{(r)}-B_{i}^{(r)}\big|\bigg)+2\sum_{r=1}^{k}
\bigg (\displaystyle\sup_{\genfrac{}{}{0pt}{}{0\leq s,t\leq
N}{|s-t|<2}}\big|B_{s}^{(r)}-B_{t}^{(r)}\big|\bigg ).\nonumber\\
\end{eqnarray}
Let us try to find an exponentially decreasing upper bound for each
term of (\ref{E12}). For this, set
$$Y_k=2\displaystyle\sum_{r=1}^{k}\bigg(\displaystyle
\max_{i=1,\ldots,N}\Big|S_{i}^{(r)}-B_{i}^{(r)}\Big|\bigg)$$
and
$$
Z_k=2\sum_{r=1}^{k}\bigg (\displaystyle\sup_{\genfrac{}{}{0pt}{}{0\leq
s,t\leq N}{|s-t|<2}}\Big|B_{s}^{(r)}-B_{t}^{(r)}\Big|\bigg ).$$ Then
\begin{eqnarray}\label{E13}
\mathbb{P}[B]&\leq& \mathbb{P} \big[Y_k+Z_k\geq 2\sqrt{N k}
(\varepsilon-\varepsilon_{1})\big]\nonumber\\ \nonumber\\
&\leq& \mathbb{P} \big [Y_k\geq
\sqrt{N k}(\varepsilon-\varepsilon_{1})\big]+\mathbb{P} \big [Z_k\geq
\sqrt{N k}(\varepsilon-\varepsilon_{1})\big].
\end{eqnarray}
 \noindent In view of the Gaussian hypothesis and the coupling, $Y_k=0$.
 Applying the Markov inequality gives for any $\lambda>0$,
\begin{eqnarray}\label{E15bis}
\mathbb{P} \big[Z_k\geq \sqrt{N k}(\varepsilon-\varepsilon_{1})\big ]&\leq&
\mathbb{E}\big[\exp{\big(\lambda Z_k^2\big)}\big].\exp{\big(-\lambda(\varepsilon-\varepsilon_1)^2N k}\big)
\nonumber\\ \nonumber\\
&\leq& \mathbb{E}\Big[\exp{\Big(2\lambda
k\Big(\displaystyle\sup_{\genfrac{}{}{0pt}{}{0\leq s,t\leq
N}{|s-t|<2}}\big|B_{s}^{(1)}-B_{t}^{(1)}\big|\Big)^2\Big)}\Big]^k
.\exp{\big(-\lambda(\varepsilon-\varepsilon_1)^2N k}\big)\nonumber\\ \nonumber\\
&\leq&\frac{\bigg(\int_0^\infty 4\lambda kt. \exp{\big(2\lambda k
t^2\big)}.\mathbb{P}\Big[\displaystyle\sup_{\genfrac{}{}{0pt}{}{0\leq
s,t\leq N}{|s-t|<2}}\big|B_{s}^{(1)}-B_{t}^{(1)}\big|\geq
t\Big]dt\bigg)^k}
{\exp{\big(\lambda(\varepsilon-\varepsilon_1)^2N k}\big)}.\nonumber \\
\end{eqnarray}
However,
\begin{equation*}
\begin{split}
\mathbb{P}\Big [\displaystyle\sup_{\genfrac{}{}{0pt}{}{0\leq
s,t\leq N}{|s-t|<2}} |B_{s}^{(1)}-B_{t}^{(1)} |\geq
t\Big ] &\leq
 \sum_{i=0}^{N-3}\mathbb{P}\Big [\sup_{i\leq t\leq
i+3}B_{t}-\inf_{i\leq t\leq i+3}B_{t}\geq
t\Big ]\\\\
&\leq N\mathbb{P}\Big [\displaystyle\sup_{0\leq t\leq
3}|B_{t}|\geq
t/2\Big ].
\end{split}
\end{equation*}
By the Brownian motion reflection principle (see for example
\cite{revuz1999cma}), $\sup_{0\leq t\leq
a}B_{t}\overset{d}{=}|B_{a}|$. Thus,
\begin{eqnarray}\label{E15}
&\mathbb{P}\Big [\displaystyle\sup_{\genfrac{}{}{0pt}{}{0\leq
s,t\leq N}{|s-t|<2}} |B_{s}^{(1)}-B_{t}^{(1)} |\geq
t\Big ]&\leq
4N\mathbb{P}\Big [B_{3}\geq
t/2\Big ]\nonumber\\ \nonumber\\
&&\leq C_{1}N\exp{\Big ( -\frac{t^2}{C_1}\Big ) }.
\end{eqnarray}
where $C_1$ is a numerical positive constant. Now, we insert
(\ref{E15}) in the integral in (\ref{E15bis}) and we choose
$\lambda=\frac{c}{k}$ were $c$ is a positive constant smaller than
$\frac{2}{C_1}$. Then we get \
\begin{equation}\label{E15biss}
\mathbb{P} \big[Z_k\geq \sqrt{N k}(\varepsilon-\varepsilon_{1})\big
]\leq C_2\exp{\Big(\frac{k\log{N}}{C_2}\Big)}.
\exp{\Big(-\frac{(\varepsilon-\varepsilon_1)^2N}{C_2}\Big)}.
\end{equation}

\

\noindent The last bound leads to the
condition $k=\mathrm{o}(N/\log{N})$. Combining (\ref{E10}), (\ref{E11}),
(\ref{E11bis}) and (\ref{E15biss}) then leads to
\begin{equation}\label{E16}
\mathbb{P}[A]\leq \mathbb{P}\big [L(1,k) \geq
2\sqrt{k}(1+\varepsilon_{1})\big]+C_2
\exp{\Big(-\frac{(\varepsilon-\varepsilon_1)^2N-k\log{N}}{C_2}\Big)}
\end{equation}
and
\begin{equation}\label{E16bis}
\mathbb{P}[A]\geq \mathbb{P}\big [L(1,k) \geq
2\sqrt{k}(1+2\varepsilon-\varepsilon_{1})\big]-C_2
\exp{\Big(-\frac{(\varepsilon-\varepsilon_1)^2N-k\log{N}}{C_2}\Big)}.
\end{equation}

\

\noindent Dividing (\ref{E16}) by
$e^{-kJ_{GUE}(\varepsilon_{1})}$ and (\ref{E16bis}) by
$e^{-kJ_{GUE}(2\varepsilon-\varepsilon_{1})}$, taking their logarithm and
then dividing the results by $k$, we get for $\alpha <\frac{1}{2}$,
\begin{equation}\label{E17}
\frac{1}{k}\log{\bigg
(\frac{\mathbb{P}[A]}{e^{-kJ_{GUE}(\varepsilon_{1})}}\bigg)}\leq
\frac{1}{k}\log{\bigg (\frac{\mathbb{P}[L(1,k) \geq
2\sqrt{k}(1+\varepsilon_{1})]}{e^{-kJ_{GUE}(\varepsilon_{1})}}+g_k(\varepsilon_1,\varepsilon)\bigg
)}
\end{equation}
and
\begin{equation}\label{E17bis}
\frac{1}{k}\log{\bigg
(\frac{\mathbb{P}[A]}{e^{-kJ_{GUE}(2\varepsilon-\varepsilon_{1})}}\bigg
)}\geq\frac{1}{k}\log{\bigg (\frac{\mathbb{P}[L(1,k)\geq
2\sqrt{k}(1+2\varepsilon-\varepsilon_{1})]}{e^{-kJ_{GUE}
(2\varepsilon-\varepsilon_{1})}}-g'_k(\varepsilon_1,\varepsilon)\bigg
)}
\end{equation}
where

\

$$
 g_k(\varepsilon_1,\varepsilon)=C_2\exp{\Big (
-\frac{(\varepsilon-\varepsilon_{1})^{2}N}{C_2}
+k\big(\frac{\log{N}}{C_2}+J_{GUE}(\varepsilon_{1})\big)\Big)}
$$
and
$$
g'_k(\varepsilon_1,\varepsilon)=C_2\exp{\Big (
-\frac{(\varepsilon-\varepsilon_{1})^{2}N}{C_2}
+k\big(\frac{\log{N}}{C_2}+J_{GUE}(2\varepsilon-\varepsilon_{1})\big)\Big)},
$$

\

 \noindent
 are two positive functions. Moreover, for $k$ large enough,
 $g_k(\varepsilon_1,\varepsilon)$ and $g'_k(\varepsilon_1,\varepsilon)$ are negligible with respect to
 $e^{-\eta k}$ for every $\eta>0$ since $k=\mathrm{o}(N/\log{N})$.
 Thus, using (\ref{E2}), a straightforward computation
 shows that the right-hand sides of (\ref{E17}) and  (\ref{E17bis}) both converge to zero
when $k\rightarrow \infty$. In other words, for $k=\mathrm{o}(N/\log{N})$
and $\varepsilon_{1}<\varepsilon$,
\begin{equation}\label{E18}
\limsup_{k\to \infty}\frac{1}{k}\log{\mathbb{P}[A]}\leq
-J_{GUE}(\varepsilon_{1})
\end{equation}
and
\begin{equation}\label{E19}
\liminf_{k\to \infty}\frac{1}{k}\log{\mathbb{P}[A]}\geq
-J_{GUE}(2\varepsilon-\varepsilon_{1}).
\end{equation}
Finally, notice that $J_{GUE}(\varepsilon)$ is a continuous function of
$\varepsilon>0$. It therefore follows from (\ref{E18}) and (\ref{E19})
that for every $\varepsilon >0$,
$$
\lim_{k\to \infty}\frac{1}{k}\log{\mathbb{P}[G(N,k)\geq
2\sqrt{N k}(1+\varepsilon)]}= -J_{GUE}(\varepsilon).
$$

\

The proof of the leftmost charge formula is similar. Set now, for
all $\varepsilon>0$ and $\varepsilon_{1}<\varepsilon$,
$$
E=\big \{G(N,k)\leq 2\sqrt{N k}(1-\varepsilon)\big \}.
$$
By the same arguments as before, we get
\begin{equation}\label{E20}
\mathbb{P}[E]\leq \mathbb{P}\big [L(N,k)\leq
2\sqrt{N k}(1-\varepsilon_{1})\big ]+\mathbb{P}[B]
\end{equation}
and
\begin{equation}\label{E21}
\mathbb{P}[E]\geq \mathbb{P}\big [L(N,k)\leq
2\sqrt{N k}(1-2\varepsilon+\varepsilon_{1})\big ]-\mathbb{P}[B].
\end{equation}
Furthermore, by (\ref{E3}),
\begin{equation}\label{E22}
\lim_{k\to \infty}\frac{1}{k^{2}}\log{\bigg
(\frac{\mathbb{P}[L(1,k) \leq
2\sqrt{k}(1-\varepsilon_{1})]}{e^{-k^{2}I_{GUE}(\varepsilon_{1})}}\bigg
)}=0
\end{equation}
and
\begin{equation}\label{E22bis}
\lim_{k\to \infty}\frac{1}{k^{2}}\log{\bigg
(\frac{\mathbb{P}[L(1,k) \leq 2\sqrt{k}(1-2\varepsilon+\varepsilon_{1})
]}{e^{-k^{2}I_{GUE}(2\varepsilon-\varepsilon_{1})}}\bigg )}=0 .
\end{equation}
\

\noindent Using the same upper bound on $\mathbb{P}[B]$ and combining
(\ref{E20}), (\ref{E21}), (\ref{E22}) and (\ref{E22bis}), one can
easily deduce that, for $k=\mathrm{o}(N^{1/2})$ and
$\varepsilon_{1}<\varepsilon$,
\begin{equation*}
\limsup_{k\to \infty}\frac{1}{k^{2}}\log{\mathbb{P}[E]}\leq
-I_{GUE}(\varepsilon_{1})
\end{equation*}
and
\begin{equation*}
\liminf_{k\to \infty}\frac{1}{k^{2}}\log{\mathbb{P}[E]}\geq
-I_{GUE}(2\varepsilon-\varepsilon_{1}).
\end{equation*}
At this stage, let us assume that $I_{GUE}(\varepsilon)$ is a
continuous function of $\varepsilon$. Then, for $k=\mathrm{o}(N^{1/2})$,
$$
\lim_{N\rightarrow \infty}\frac{1}{k^{2}}~\log{\mathbb{P}\big
[G(N,k) \leq 2\sqrt{N k}(1-\varepsilon)\big ]}= -I_{GUE}(\varepsilon),
$$
which is the result.


We are left with the proof of the continuity of $I_{GUE}(\varepsilon)$.
Set $\mathcal{M}((-\infty,t])$, the set of all probability measures
on $(-\infty,t]$ when $ t \in \mathbb{R}$. For a given distribution
$\mu\in\mathcal{M}((-\infty,t])$, define the corresponding potential
energy, as in \cite{saff}, by
\begin{equation*}
I_{\mu}(t)=2\int_{-\infty}^{t}x^{2}d\mu(x)-\int_{-\infty}^{t}\int_{-\infty}^{t}\log{|x-y|}d\mu(x)d\mu(y).
\end{equation*}
The minimal energy
$$
I(t)=\inf_{\mu\in\mathcal{M}((-\infty,t])}I_{\mu}(t)
$$
precisely allows us to compute the rate function $I_{GUE}$ via the
formula $I_{GUE}(\varepsilon)=I(1-\varepsilon)-I(\infty)$. The last
equality could be found in \cite{feral2004lds}. For $t\geq 1$,
$I(t)$ is a constant function, the extremal measure is the so-called
semi-circular law supported on $[-1,1]$  and the energy
$I(t)=\log{(2)}+3/4$, (cf. \cite{benalic,saff}). For each
$t\in\mathbb{R}$, there is a unique measure
$\nu_{t}\in\mathcal{M}((-\infty,t])$, with no mass point, achieving
the infimum (cf. \cite{saff}). Furthermore, $\nu_{t}$ is compactly
supported and the corresponding energy is finite. Since $I(t)$ is an
infimum and a non-increasing function of $t$, for any $\eta>0$,
\begin{equation}\label{C2}
 I(t)\leq
I(t-\eta)\leq
\frac{I_{\nu_{t}}(t-\eta)}{\nu^{2}_{t}((-\infty,t-\eta])} \, .
\end{equation}
It is obvious that the right-hand side of (\ref{C2}) converges to
$I(t)$ when $\eta$ converges to zero. This proves the
left-continuity of $I(t)$.

To show the right-continuity, notice that by a simple change of
variable,
\begin{equation*}
I(t)=\inf_{\mu\in\mathcal{M}((-\infty,t+\eta])}I^{\eta}_{\mu}(t)
\end{equation*}
where
\begin{equation*}
I^{\eta}_{\mu}(t)=2\int_{-\infty}^{t+\eta}(x-\eta)^{2}d\mu(x)-\int_{-\infty}^{t+\eta}
\int_{-\infty}^{t+\eta}\log{|x-y|}d\mu(x)d\mu(y).
\end{equation*}
 Consequently,
\begin{eqnarray}\label{C4}
 &I(t)-I(t+\eta)&\leq I^{\eta}_{\nu_{t+\eta}}(t)-I(t+\eta)\nonumber \\
 \nonumber \\
&&\leq 2\eta^{2}+4\eta\int_{-\infty}^{t+\eta}|x|~d\nu_{t+\eta}(x).
\end{eqnarray}
For $|x|\geq |y|$, we have $\log{|x-y|}\leq \log{|2x|}$. Moreover,
there is a positive constant $C_{3}$ such that $|x|\leq
C_{3}(2x^{2}-2\log{|2x|})$. In view of (\ref{C4}), a straightforward
calculation leads to
\begin{equation}\label{C5}
I(t)-I(t+\eta)\leq 2\eta^{2}+4\eta ~C_{3} I(t).
\end{equation}
Since $I(t)$ is finite, the right-hand side of (\ref{C5}) converges
to zero when $\eta \rightarrow 0$. Thus, the continuity of $I(t)$ is
proved, and  that of $I_{GUE}(\varepsilon)$ as well. The proof of
Theorem \ref{th1} is now complete. \vspace{0.4cm}

\hspace{16cm}$\square$


\section{Exponential-tailed distribution and the KMT approximation}\label{sec3}
\setcounter{equation}{0} In this section, we replace the  standard
normal variables with weights having finite moment-generating
function around zero. When comparing $G(N,k)$ to the Brownian
last-passage percolation model, $Y_k$ will not vanish as in the
Gaussian case where the coupling was "perfect". Actually, we couple
easily a partial sum of $i.i.d.$ standard normal random weights with
a Browian motion. However, loosing this assumption will complicate
the task and a new coupling is then required. Following
\cite{bodineau2005upl}, we make use of the KMT approximation: a
powerful tool to couple a partial sum of $i.i.d.$ random variables
and a Wiener process, both constructed on the same probability
space. The KMT approximation, also called the invariance principle,
was first introduced in 1975 by Koml\'os, Major and Tusn\'ady in the
famous work \cite{komlos}. The basic version deals with a partial
sum of $i.i.d.$ random variables reconstructed in a way to be
''close'' to another partial sum of $i.i.d.$ standard normal random
variables. Later versions of this strong approximation do not
require a common distribution, see \cite{sakhanenko}. The readers can
also see \cite{csorgo} for a complete survey.

Let $(X_i)_{i\geq1}$ be a sequence of independent random variables
and denote by $S_N$ the corresponding partial sum. Let
$(B_t)_{t\geq0}$ be a Brownian motion built on the same probability
space. The following theorem is an immediate consequence of Theorem
1 in \cite{komlos}.

\begin{theorem}[Koml\'os-Major-Tusn\'ady]\label{th10}
Assume that $\mathbb{E}X_1=0$ and $\mathbb{E}X_1^{2}=1$. Assume
further that there exit $\mu_0>0$ such that for all $\mu<\mu_0$,
\begin{equation*}
\mathbb{E}\exp{(\mu|X_1|)}<+\infty.
\end{equation*}
Then for every $N\geq 1$,  the sequence $(X_i)_{i\geq1}$ and the
Brownian motion $(B_t)_{t\geq0}$ can be constructed in such a way
that for all $x>0$,
$$
\mathbb{P}\big[\max_{i=1,\ldots,N}\big|S_i-B_i\big|>\Theta\log{N}+x\big]\leq
C\exp{(-\theta x)}.
$$
The positive constants $\Theta$, $C$ and $\theta$ depend only on the
distribution of $X_1$ and $\theta$ can be taken as large as desired
by choosing $\Theta$ large enough.
\end{theorem}

According to the notation and the steps of Section \ref{sec2},
recall that
\begin{equation*}
|G(N,k)-L(N,k)|\leq Y_k+Z_k .
\end{equation*}
We already have $Z_k$ in (\ref{E15biss}) and we want $Y_k$ to be as
small as possible. To this, we construct the sequence
$(X_{i}^{(j)})_{i,j=1}^{\infty}$ and the independent Brownian
motions $(B_t^{(r)})_{t\geq0}$ in the sense of
 Theorem \ref{th10}. By the Markov inequality, we have for all $\varepsilon>0$,
 $\varepsilon_1<\varepsilon$ and $\lambda>0$,

\begin{eqnarray*}
\mathbb{P} \big[Y_k\geq \sqrt{N k}(\varepsilon-\varepsilon_{1})\big
]&\leq& \mathbb{E}\big[\exp{\big(\lambda
Y_k\big)}\big].\exp{\big(-\lambda(\varepsilon-\varepsilon_1)\sqrt{N k}}\big)
\nonumber\\ \nonumber\\
&\leq& \mathbb{E}\Big[\exp{\Big(2\lambda
\max_{i=1,\ldots,N}\Big|S_{i}^{(1)}-B_{i}^{(1)}\Big|\Big)}\Big]^k
.\exp{\big(-\lambda(\varepsilon-\varepsilon_1)\sqrt{N k}}\big)\nonumber\\ \nonumber\\
&\leq&\frac{\bigg(\int_0^\infty 2\lambda . \exp{\big(2\lambda
t\big)}.\mathbb{P}\Big[\max_{i=1,\ldots,N}\Big|S_{i}^{(1)}-B_{i}^{(1)}\Big|\geq
t\Big]dt\bigg)^k}
{\exp{\big(\lambda(\varepsilon-\varepsilon_1)\sqrt{N k}}\big)}.\nonumber \\
\end{eqnarray*}

In order to apply Theorem \ref{th10}, we make the simple variable
change $t=s-\Theta\log{N}$ and we choose $\lambda<\theta/2$.
Therefore, there exist two positive constant $c_4$ and $C_4$ such
that,
\begin{equation}\label{e1}
\mathbb{P} \big[Y_k\geq \sqrt{N k}(\varepsilon-\varepsilon_{1})\big
]\leq
C_4\exp{\Big(-\frac{(\varepsilon-\varepsilon_1)\sqrt{N k}-k\log{N}}{C_4}\Big)}.
\end{equation}
Now, putting  (\ref{E10}), (\ref{E11}), (\ref{E15biss}) and
(\ref{e1}) together gives

\begin{equation}\label{e2}
\begin{split}
\mathbb{P}[A]\leq &\mathbb{P}\Big[L(N,k) \geq
2\sqrt{N k}(1+\varepsilon_{1})\Big]+C_5
\exp{\Big(-\frac{(\varepsilon-\varepsilon_1)^2N-k\log{N}}{C_5}\Big)}\\
&+C_5\exp{\Big(-\frac{(\varepsilon-\varepsilon_1)\sqrt{N k}-k\log{N}}{C_5}\Big)}
\end{split}
\end{equation}
and
\begin{equation}\label{e3}
\begin{split}
\mathbb{P}[A]\geq &\mathbb{P}\Big[L(N,k) \geq
2\sqrt{N k}(1+2\varepsilon -\varepsilon_{1})\Big]-C_5
\exp{\Big(-\frac{(\varepsilon-\varepsilon_1)^2N-k\log{N}}{C_5}\Big)}\\
&-C_5\exp{\Big(-\frac{(\varepsilon-\varepsilon_1)\sqrt{N k}-k\log{N}}{C_5}\Big)}.
\end{split}
\end{equation}
On the left of the mean, we have
\begin{equation}\label{e4}
\begin{split}
\mathbb{P}[E]\leq &\mathbb{P}\Big[L(N,k) \leq
2\sqrt{N k}(1-\varepsilon_{1})\Big]+C_5
\exp{\Big(-\frac{(\varepsilon-\varepsilon_1)^2N-k\log{N}}{C_5}\Big)}\\
&+C_5\exp{\Big(-\frac{(\varepsilon-\varepsilon_1)\sqrt{N k}-k\log{N}}{C_5}\Big)}
\end{split}
\end{equation}
and
\begin{equation}\label{e5}
\begin{split}
\mathbb{P}[E]\geq &\mathbb{P}\Big[L(N,k) \leq
2\sqrt{N k}(1-2\varepsilon+\varepsilon_{1})\Big]-C_5
\exp{\Big(-\frac{(\varepsilon-\varepsilon_1)^2N-k\log{N}}{C_5}\Big)}\\
&-C_5\exp{\Big(-\frac{(\varepsilon-\varepsilon_1)\sqrt{N k}-k\log{N}}{C_5}\Big)}.
\end{split}
\end{equation}
Proceeding like in Section \ref{sec2}, we divide (\ref{e2}) by
$e^{-kJ_{GUE}(\varepsilon_{1})}$, (\ref{e3}) by
$e^{-kJ_{GUE}(2\varepsilon-\varepsilon_{1})}$, (\ref{e4}) by
$e^{-k^2I_{GUE}(\varepsilon_{1})}$ and (\ref{e5}) by
$e^{-k^2I_{GUE}(2\varepsilon-\varepsilon_{1})}$. To handle the
remaining parts when $k\to \infty$, we take $k$ negligible with
respect to $N$. On the right of the mean, we need
$k=\mathrm{o}\big(\frac{N}{(\log{N})^2}\big)$ and on the left, we
take $k=\mathrm{o}(N^{1/3})$. Finally we conclude as in section
\ref{sec2} using the continuity of $J_{GUE}$ and $I_{GUE}$.
 \vspace{0.4cm}

\hspace{16cm}$\square$


\section{Subexponential weights and the Skorokhod embedding}\label{sec4}
\setcounter{equation}{0} In this section, we consider a particular
category of subexponential weights verifying $\mathbb{E}\exp{(\mu|X|^\gamma)}\\<\infty$
for some $\mu>0$ and $\gamma\in(0,1)$. We say that these variables have a Weibull-like
tails because of the similarity with the right tail of the Weibull distribution with
a shape parameter lower than $1$. Such weights are considered to be heavy-tailed
and then, do not satisfy the Cram\'er condition. However, $G(N,k)$ still satisfies an LDP
principle on a very thin rectangle. For subexponential weights with finite $p$-th moment for $p\geq2$,
we need $k$ to be at least smaller than $\log{(N^{p-1})}$ in order to prove deviation results.
Computations in this case are left to the reader. The precise definition of a subexponential distribution and
large deviations for a partial sum of $i.i.d.$ such weights can be found in \cite{nagaev}.

We do not know a KMT strong approximation version for random variables satisfying the moment condition above.
We follow \cite{baik2005gcl} by using the Skorokhod embedding theorem instead. The Skorokhod embedding theorem
is another tool to couple a sum of $i.i.d.$ random variables with a Brownian motion, see (cf.
\cite{breiman, skorokhod, obloj}) for more details concerning this
theorem.


\begin{theorem}[Skorokhod]\label{th3}
Let $(B_{t})_{t\geq 0}$ be a standard one-dimensional Brownian motion
and $X$ a real valued random variable satisfying $\mathbb{E}X=0$ and
$\mathbb{E}X^{2}=1$. Then, there is a stopping time $T$ for the
Brownian motion such that $B_{T}\overset{d}{=}X$ and $\mathbb{E}T=1$.
\end{theorem}

An immediate consequence of this theorem allows to embed sums of
real independent random variables into the Brownian motion. Applying
the strong Markov property to the Brownian motion, Theorem \ref{th3}
yields the following classical corollary.

\begin{corollary}\label{cor2}
Let $X_{1},X_{2},\ldots,X_{N},\ldots$ be i.i.d. satisfying
$\mathbb{E}X_{1}=0$, $\mathbb{E}X_{1}^{2}=1$ and set
$S_{N}=X_{1}+X_{2}+\cdots+X_{N}$, $N \geq 1$. There is a sequence of
 i.i.d. stopping times $\tau_{0}=0,\tau_{1},\ldots,\tau_{N},\ldots$
such that
$$
S_{N}\overset{d}{=}B_{\tau_{1}+\cdots+\tau_{N}}
$$
and $
(B_{\tau_{1}+\cdots+\tau_{N+1}}-B_{\tau_{1}+\cdots+\tau_{N}})_{N
\geq 0} $ is a sequence of i.i.d. random variables having the same
distribution as $X_{1}$.
\end{corollary}

In our context, an application of the last corollary allows to claim
that there exists $i.i.d.$ stopping times for the Brownian motion
$\tau_{0}=0,\tau_{1},\ldots,\tau_{N},\ldots$ such that
$\mathbb{E}\tau_{1}=1$ and
$S_{i}^{(r)}\overset{d}{=}B^{(r)}_{\tau_{1}+\cdots+\tau_{i}}$ for
$i\geq 1$ and $r\geq 1$. Consequently, choose
$$
X_{i}^{(r)}=B^{(r)}_{\tau_{1}+\cdots+\tau_{i}}-B^{(r)}_{\tau_{1}+\cdots+\tau_{i-1}},
$$
in order to have
$$
S^{(r)}_{i}=B^{(r)}_{\tau_{1}+\cdots+\tau_{i}}\,\,\,\,\,\,\,\,\,\,\,a.s.
$$
Thus, for any $t\geq0$,
\begin{equation*}
\begin{split}
&\mathbb{P}\Big [\displaystyle\max_{i=1,\ldots,N}\big
|S^{(1)}_{i}-B_{i}^{(1)}\big|\geq t\Big]
\\\\
= &\mathbb{P}\Big
[\displaystyle\max_{i=1,\ldots,N}\big
|B^{(1)}_{\tau_{1}+\ldots+\tau_{i}}-B_{i}^{(1)}\big |
\geq t;
\displaystyle\max_{i=1,\ldots,N}\Big |\sum_{l=1}^{i}(\tau_{l}-1)\Big
|\geq
t^{\beta}\Big ]\\
\quad &+\mathbb{P}\Big
[\displaystyle\max_{i=1,\ldots,N}\big
|B^{(1)}_{\tau_{1}+\ldots+\tau_{i}}-B_{i}^{(1)}\big |
\geq t;
\displaystyle\max_{i=1,\ldots,N}\Big |\sum_{l=1}^{i}(\tau_{l}-1)\Big
|< t^{\beta}\Big ].\\\\
\end{split}
\end{equation*}
Hence
\begin{eqnarray}\label{E24}
&\mathbb{P}\Big [\displaystyle\max_{i=1,\ldots,N}\big
|S^{(1)}_{i}-B_{i}^{(1)}\big|\geq t\Big] &\leq
\mathbb{P}\Big[\displaystyle\sup_{\genfrac{}{}{0pt}{}{0\leq
s,t\leq N}{|s-t|<N^{\beta}}}\Big|B_{s}^{(1)}-B_{t}^{(1)}\Big|\geq
t\Big]\nonumber\\
&&\,\,\,+\mathbb{P}\Big[\max_{i=1,\ldots,N}\Big|\sum_{l=1}^{i}(\tau_{l}-1)\Big|\geq
t^{\beta}\Big]. \\\nonumber
\end{eqnarray}
We evaluate each term of (\ref{E24}) separately. First,
\begin{eqnarray*}
&&\mathbb{P}\Big[\displaystyle\sup_{\genfrac{}{}{0pt}{}{0\leq
s,t\leq N}{|s-t|<t^{\beta}}}\Big|B_{s}^{(1)}-B_{t}^{(1)}\Big|\geq
t\Big]\\ \\\\
&&\leq\sum_{i=0}^{N-t^{\beta}}\mathbb{P}\Big[\sup_{i\leq t\leq
i+t^{\beta}+1}B_{t}-\inf_{i\leq t\leq i+t^{\beta}+1}B_{t}\geq
t\Big]\\ \\\\
 &&\leq
N\mathbb{P}\Big[\sup_{0\leq t\leq t^{\beta}+1}|B_{t}|\geq
t/2\Big].
\end{eqnarray*}
Applying the reflection principle as in Section \ref{sec2}, we get
\begin{eqnarray}\label{E25}
&\mathbb{P}\Big[\displaystyle\sup_{\genfrac{}{}{0pt}{}{0\leq
s,t\leq N}{|s-t|<t^{\beta}}}\Big|B_{s}^{(1)}-B_{t}^{(1)}\Big|\geq
t\Big] &\leq
4N\mathbb{P}\Big[B_{t^{\beta}+1}\geq
t/2\Big]\nonumber\\
\nonumber\\
\nonumber\\&&\leq4N\exp{\Big(-\frac{t^{2-\beta}}{8}\Big)}.
\end{eqnarray}
To find an upper bound for the second term on the right-hand side of (\ref{E24}),
we need a connection between the weight moments and
those of the stopping times obtained by the Skorokhod embedding.
Furthermore, we need to control the sum of the independent stopping
times to reach an exponentially decaying inequality. When the weights are
bounded for exemple, we can construct a stopping time with finite exponential moments.
The sum is then controlled by the Bernstein inequality.

However, when $X$ only satisfies (\ref{E46}), the Skorokhod stopping
time does not necessarily have a finite exponential moment
and thus the Bernstein inequality can not be applied. For example,
in \cite{davis}, Davis found the best universal constant connecting
the stopping time moments to those of the stopped Brownian motion.
More precisely, if $(B_{t})_{t>0}$ is a Brownian motion and $\tau$
is a stopping time, then there is a universal constant $a_{p}$ such
that, when $1<p<\infty$ and $\mathbb{E}\tau^{p/2}<+\infty$,
\begin{equation}\label{E50}
a_{p}\mathbb{E}\tau^{p/2}\leq\mathbb{E}|B_{\tau}|^{p}.
\end{equation}
 Moreover, Davis proved that
the best constant for $p=2n$ ($n\in\mathbb{N}^{*}$) is
${z_{2n}^{*}}^{2n}$ which is the smallest positive zero of the
Hermite polynomial of order $2n$. In \cite{burkholder}, this
constant is shown to be $O((2n)^{-n}))$. So unless $X=B_{\tau}$ is a
bounded variable, $\tau$ can not have finite exponential moments.

The constant above is universal but it could be sharpened for some
particular stopping times. For example, considering the stopping
time of the Skorokhod representation \cite{breiman, skorokhod},
Sawyer improved the constant $a_{p}$ and established, in
\cite{sawyer}, the following inequality.

\begin{theorem}\label{th8}
Let $X$ be a centered random variable such that
\begin{equation*}
\mathbb{E}\exp{(\mu|X|^{\gamma})}<+\infty
\end{equation*}
  for some $\gamma>0$
and $\mu>0$, and let $\tau$ be the corresponding stopping time of
the Skorokhod representation. Set $\theta=\frac{\gamma}{2+\gamma}$ and
$\nu=\mu^{1-\theta}$. Then,
\begin{equation*}
\mathbb{E}\exp{(\nu\tau^{\theta})}\leq \Phi_{\gamma}\,
\mathbb{E}\exp{(\mu |X|^{\gamma})},
\end{equation*}
for some positive constant $\Phi_{\gamma}$ depending only on
$\gamma$.
\end{theorem}
Note that a similar exponential bound may be obtained from
(\ref{E50}). However the cost is a worse constant $\mu$.

\noindent Under the assumption (\ref{E46}) and in view of Theorem \ref{th8},
the Bernstein inequality can not be
applied to the sum of the independent stopping times because
$\theta=\gamma/(2+\gamma)<1$. To avoid this obstacle, we introduce
the Fuk-Nagaev inequality \cite{Fuk} which requires less restrictive
assumptions.
\begin{theorem}[Fuk-Nagaev]\label{th9}
Let $X_1,\cdots,X_N$ be a sequence of real i.i.d. random variables
satisfying $\mathbb{E}X_1=0$ and $\mathbb{E}X_{1}^{2}=\sigma^{2}$.
Then, for all $x>0$ and $y>0$,
\begin{equation}\label{E49}
\mathbb{P}\Big[\max_{1\leq i\leq
N}{\Big|\sum^{i}_{l=1}X_{l}\Big|}\geq x\Big]\leq
N\mathbb{P}\Big[|X_1|> y
\Big]+2\exp{\Big(-\frac{x^{2}}{2(N\sigma^{2}+xy/3)}\Big)}.
\end{equation}
\end{theorem}
We refer to \cite{fuk1} for more detail on this inequality.
Recall now the second term of the right hand side of (\ref{E24}) and
consider the stopping times of the Skorokhod representation.
Choosing $x=t^{\beta}$ and $y=t^{\delta}$ in Theorem \ref{th9} and
applying Markov inequality to the first term of the right hand side
of (\ref{E49}), one has
\begin{eqnarray}\label{E47}
&\mathbb{P}\Big[\displaystyle\max_{i=1,\ldots,N}\Big|\displaystyle\sum_{l=1}^{i}(\tau_{l}-1)\Big|\geq
t^{\beta}\Big] &\leq
C_6N\exp{\Big(-\frac{\nu
t^{\theta\delta}}{C_6}\Big)}+C_6\exp{\Big(-\frac{t^{2\beta}}{C_6\max{\{N,t^{\beta+\delta}\}}}\Big)}.
\end{eqnarray}

\noindent Now, we apply Markov inequality to $Y_k$ as in Section \ref{sec3}
and we get for some $\lambda>0$ and$0<\eta<1$,
\begin{eqnarray}\label{E48}
\mathbb{P} \big[Y_k\geq \sqrt{N k}(\varepsilon-\varepsilon_{1})\big
]&\leq& \mathbb{E}\big[\exp{\big(\lambda
Y_k^\eta\big)}\big].\exp{\big(-\lambda(\varepsilon-\varepsilon_1)^\eta (N k)^{\eta/2}}\big)
\nonumber\\  \nonumber\\
&\leq&\frac{\bigg(\int_1^\infty 4\lambda \exp{\big(2\lambda
t^\eta\big)}\mathbb{P}\Big[\max_{i=1,\ldots,N}\Big|S_{i}^{(1)}-B_{i}^{(1)}\Big|\geq
t\Big]dt\bigg)^k}
{\exp{\big(\lambda(\varepsilon-\varepsilon_1)^\eta (N k)^{\eta/2}}\big)}.\nonumber \\
\end{eqnarray}
Inserting (\ref{E24}) and (\ref{E47}) in (\ref{E48}) and choosing
$\lambda$ very small yield the following constraints on $\eta$,
$\beta$, $\delta$ and $\theta$.
$$\Bigg\{
\begin{array}{l}
 \eta<2-\beta\\
\eta<\beta-\delta\\
\eta<\theta\delta.
\end{array}
$$
 Straightforward computations lead us to choose
$\beta=\frac{4\gamma+4}{3\gamma+2}$ and
$\delta=\frac{2\gamma+4}{3\gamma+2}$ since
$\theta=\frac{\gamma}{\gamma+2}$. Consequently, we obtain
$\eta<\frac{2\gamma}{3\gamma+2}$. To get large deviation asymptotic
formulas on the right and the left of the mean, we respectively need
$\alpha<\frac{\eta}{2-\eta}$ and $\alpha<\frac{\eta}{4-\eta}$. This
completes the proof of Theorem \ref{th7}.

\vspace{0.4cm}

\hspace{16cm}$\square$


\section{Small and large deviations inequalities}\label{sec5}
\setcounter{equation}{0} Non-asymptotic bound on the right and the
left of the mean is an immediate consequence of the corresponding
bound for the GUE and the arguments developed in Sections 2 and 3.
In particular, we use (cf.\cite{ledoux}) that there exists a
positive constant $C_7$ such that, for any $\varepsilon>0$,
\begin{equation}\label{E3bis}
\mathbb{P}\big[L(1,k) \geq 2\sqrt{k}(1+\varepsilon)\big]\leq
\exp{\big(-kJ_{GUE}(\varepsilon)\big)} \leq C_7 \,
\exp{\Big(-\frac{k\max{(\varepsilon^{2},\varepsilon^{3/2})}}{C_7}\Big)}.
\end{equation}
On the left of the mean, deviation inequalities for the largest
eigenvalue of the GUE for a given $k$ are quite more complicated to
prove. Ledoux and Rider obtained in a recent paper,
\cite{ledouxrider}, that the leftmost charge of the largest
eigenvalue of a large set of random matrices behaves like the left
tail of the corresponding Tracy-Widom law. More precisely, they get
for all $0<\varepsilon\leq 1$,
\begin{equation}\label{E3biss}
\mathbb{P}\big[L(1,k) \leq 2\sqrt{k}(1-\varepsilon)\big]\leq C_7 \,
\exp{\bigg(-\frac{k^2\varepsilon^{3}}{C_7}\bigg)}.
\end{equation}

 \noindent As we mentioned before, when $\varepsilon>1$, we have
 Gaussian behavior for both left and right
tails. This follows from concentration arguments dealing with
Lipschitz functions of independent standard normal variables. Once
more, two cases will be tackled: Standard normal weights and finite
exponential moments ones.
\subsection{Standard normal variables}
Following the proof of Theorem \ref{th1} in Section \ref{sec2},
choose $\varepsilon_{1}=\frac{\varepsilon}{2}$. Then, combining
(\ref{E16}) and (\ref{E3bis}), for any $\varepsilon>0$,
\begin{equation*}
\mathbb{P}[A]\leq
C_8\exp{\Big(-\frac{k\max{(\varepsilon^{3/2},\varepsilon^{2})}}{C_8}\Big)}+C_8
\exp{\Big(-\frac{\varepsilon^{2}N-k\log{N}}{C_8}\Big)}
\end{equation*}
where $C_8 >0$. In order to reach (\ref{E8}) when $0<\varepsilon<1$,
we need a positive constant $C (\alpha)>C_8$, depending only on
$\alpha$, such that
\begin{equation}\label{E42}
C_\alpha\exp{\Big(-\frac{k\varepsilon^{3/2}}{C_\alpha}\Big)}\geq C_8
\exp{\Big(-\frac{\varepsilon^{2}N-k\log{N}}{C_8}\Big)}.
\end{equation}
Taking the logarithm of (\ref{E42}), $C_\alpha$ has to satisfy
\begin{equation}\label{E42bis}
\log{\frac{C_8}{C_\alpha}}-k\varepsilon^{3/2}
\Big(\frac{N\varepsilon^{1/2}}{C_8k}-
\frac{\log{N}}{C_8\varepsilon^{3/2}}-\frac{1}{C_\alpha}\Big)\leq0.
\end{equation}
However, since $\mathbb{P}[A]\leq 1$, $\varepsilon$ has to satisfy
\begin{equation}\label{E43}
k\varepsilon^{3/2}\geq 1.
\end{equation}
Combining now (\ref{E42bis}) and (\ref{E43}), we  finally get that
$C_\alpha$ has to satisfy
\begin{equation}\label{E44}
\log{\frac{C_8}{C_\alpha}}+\frac{1}{C_\alpha}-
\frac{N^{1-\frac{4\alpha}{3}}-N^{\alpha}\log{N}}{C_8}\leq0.
\end{equation}
Hence $C_\alpha$ exists and satisfies (\ref{E44}) only if
$\alpha<\frac{3}{7}$. In that case, we make the reverse computation
to conclude that
$$
\mathbb{P}\big[G(N,k) \geq 2\sqrt{N k}(1+\varepsilon)\big]\leq
2C_\alpha\exp{\Big(\frac{-k\varepsilon^{\frac{3}{2}}}{C_\alpha}\Big)}.
$$

\

\noindent We make the same computations for the left-tail upper
bound. Here, $C_\alpha$ has to satisfy
$$
\log{\frac{C_8}{C_\alpha}}-k^2\varepsilon^{2} \Big(\frac{N}{C_8k^2}-
\frac{\log{N}}{C_8k\varepsilon^{2}}-\frac{\varepsilon}{C_\alpha}\Big)\leq0,
$$
which finally gives
$$
\log{\frac{C_8}{C_\alpha}}+\frac{1}{C_\alpha}-\frac{N^{1-2\alpha}-N^{\alpha/3}\log{N}}{C_8}\leq0.
$$
 This proves Theorem \ref{th5}.

\hspace{16cm}$\square$
\subsection{Finite moment-generating function case}
Choosing $\varepsilon_{1}=\frac{\varepsilon}{2}$  in (\ref{e2}) and
taking into consideration (\ref{E43}), the inequalities
(\ref{E3bis}) and (\ref{E3biss}) imply that there exists a positive
constant $C_9$ depending on $\alpha$ and the distribution of $X$
such that, for all $\varepsilon> 0$,
\begin{equation*}
\begin{split}
\mathbb{P}[A]\leq
&C_9\exp{\Big(-\frac{k\varepsilon^{3/2}}{C_9}\Big)}\bigg(1+
\exp{\Big(-\frac{N^{\frac{1}{2}-\frac{\alpha}{6}}-N^\alpha\log{N}-N^\alpha}{C_9}
\Big)}\\
&\,\,\,\,\,\,\,\,\,\,\,\,\,\,\,\,\,\,\,\,\,\,\,\,\,\,\,\,\,\,\,\,\,\,\,\,\,\,\,\,\,\,\,\,
\,\,\,\,\,\,\,\, +
\exp{\Big(-\frac{N^{1-4\alpha/3}-N^\alpha\log{N}-N^\alpha}{C_9}\Big)}\bigg)
\end{split}
\end{equation*}
and
\begin{equation*}
\begin{split}
\mathbb{P}[E]\leq
&C_9\exp{\Big(-\frac{k^2\varepsilon^3}{C_9}\Big)}\bigg(1+
\exp{\Big(-\frac{N^{\frac{1-3\alpha}{2}}-1}{C_9} \Big)} +
\exp{\Big(-\frac{N^{1-2\alpha}-N^{\frac{\alpha}{3}}-1}{C_9}\Big)}\bigg).
\end{split}
\end{equation*}
This means that we have a right-tail bound for $\alpha<3/7$ and a
left-tail bound for $\alpha<1/3$. The proof is complete and thus
Theorem \ref{th6} is proved.

\hspace{16cm}$\square$

\section*{Acknowledgement} I would like to thank Michel Ledoux for many valuable
discussions on random growth and random matrices. I also thank Delphine
F\'eral for remarks on the rate function of the GUE ensemble.\\\\\\\\
\noindent \textit{Jean-Paul Ibrahim.}\\
\noindent \textit{Institut de Mathématiques, Université Paul-Sabatier, $31062$
Toulouse, France.}\\
E-mail: \verb"jibrahim@math.univ-toulouse.fr"

\bibliographystyle{abbrv}
\bibliography{LD}
\end{document}